\documentclass[12pt]{article}
\usepackage[cp866]{inputenc}
\usepackage{amsmath}
\begin{document}
\title {Prime numbers of a of a kind $n^2+1$ }
\author {V.E. Govorov}
\maketitle

\medskip
\qquad {\it Comment: 4 pages}

\medskip
\qquad {\it Subjects\,  Number Theory.}

\medskip
\qquad {\it Summary.}\,The set of prime integers of a kind $n^2+1$ is infinite.

\medskip
\qquad{\it ACM-class\, 11N05 }

\bigskip

 Let's consider a sequence  $S=$

 \begin{tabular}{cccc}
   \{
   5,                 & 17,            & 37,             & $65=5\cdot 13$,\\
   101,               &$145=5\cdot 29$,&197,             & 257,\\
   $325=5^2\cdot 13,$
    & 401,           &$485=5\cdot 97$, & 577,\\
   677,               &$785=5\cdot157$,&$901=17\cdot53$, &$1025=5^2\cdot41$\\
   $1157=13\cdot89$,  &1297,           & $1445=5\cdot17^2$,&1601,\\
   $1765=5\cdot 353$, &$1937=13\cdot149,$
   &$2117=29\cdot73$, &$2305=5\cdot461$,\\
   $2501=41\cdot61$,  &$2705=5\cdot 541$, & 2917,          & 3137, \\
   $3365=5\cdot 673$,  &$3601=13\cdot277$,\, \dots \}
 \end{tabular}
 of odd numbers kind $n^2+1$,  i.e. the numbers kind $4n^2+1$.

In additions to $S$ we'll investigate the set $P$ of prime numbers of a kind
$4k+1$ and the set  $PP$ of all numbers being products of numbers from $P$.
\bigskip

{\it Lemma 1.}\,\label{m1}{\it The set of divisors of sequence $S$ is equal to
$PP$.}

 \medskip
 {\it Proof.\,}
 Let  $m=p_1^{\alpha_1} p_2^{\alpha_2}\dots p_k^{\alpha_k}$ be the prime factor
  decomposition of  number $m$ then comparison $x^2\equiv -1 (m)$ solved  in that
   and only that case , when
 $$(-1)^{(p_i-1)/2}\equiv 1(p_i)$$ ([1], proposition 5.1.1), i.e. when
$p_i=4l+1$.

\bigskip
 {\it Lemma  2.}\,\label{m2}{\it If $S_n$ is mutuality disjoint with all
preceding numbers in the sequence $S$, then it is a prime number.}

 \medskip
{\it Proof.} Let $S_n=4n^2+1$ is mutuality disjoint with all preceding numbers
in sequence  $S$ and $d|S_n$. The number  $d$ is odd, not equal  $n$ and it may
be chosen   $d<2n$.
 if $3\leq d\leq n-1$, then we suppose $k=n-d$, if
 $n+1\leq d\leq 2n-1$ then $k=d+n$. After substituting  $k$ in
    identity :
 $$ 4n^2+1-4(n-k)(n+k)=4k^2+1$$
 we have $d|4k^2+1$,  contradicting with initial supposition.

\bigskip
  {\it Lemma 3. }\,\label{m3} {\it If $p|S_n$, then $p|S_{n+kp}$ and $p|S_{kp-n}$
  for all  $k$. If $p$ is a prime number then it divide only them.}

\medskip
 {\it Proof.}\, Let be  $4n^2+1=ps$ then
 $S_{n+kp}=4(n+kp)^2+1=4n^2+1+8nkp+4k^2p^2=ps+p(8nk+4k^2p)$.
 also $S_{kp-n}=ps+p(-8nk+4k^2p^2)$.

  Let be $p$ -- prime,  $p|S_n$ and $p|S_{n+q}=4n^2+1+8nq+4q^2=
  ps+4q(2n+q)$, then $q=k_1p$ or $2n+q=k_2p$. The first case is
  $n+q=k_1p+n$, the second case is $n+q=k_2p-n$.

\bigskip
 {\it Definition.}\,\label{d1}
 Let's define two functions    $r(m)$ and $x(m)$  on $PP$.  Function $r(m)$
 is equal to the least natural solution of  a congruence modulo  $m$  $4z^2+1 \equiv 0\,(m)$,
 i.e. $r(m)$ is a first entry number $m$ in $S$ as a factor.
  $$x(m)=\frac{4r^2(m)+1}{m}$$

   If $p=4k+1$ is a prime number  then $r(p)$ can be defined as follow. Let $g$ be a primitive
   root  modulo $p$, If   $t\equiv g^k \,(p)$ is even then  $r(p)=t/2$ else
   $r(q)=(p-t)/2$.

\bigskip
{\it Remark.}\,\label{m4} Let be  $q=p_1p_2\dots p_k$, where  $p_i$ is a prime
number of kind $4l+1$. The Ring $Z_q$  is a direct sum of fields $Z_p$ and
every fields has two square roots from  $-1$, i.e. $Q$ has $2^k$ roots from
$-1$. In this case we'll consider  $2^k$ functions $r_i(q)$ fulfilling the
congruence  $4r_i(q)^2+1\equiv 0(q)$ .

\bigskip
{\it Lemma 4.}\,\label{m5}{\it Functions  $r(m)$ and $x(m)$ satisfy to
inequations}
$$
\frac{\sqrt{m-1}}{2}\leq r(m)\leq\frac{m-1}{2},\quad x(m)\leq m-2
$$

\medskip
{\it Proof.} Since $r^2(m) \equiv (m-r(m))^2\,(m)$, and $r(m)$ is the least
integer satisfying to congruence  $4z^2\equiv -1(m)$, so $r(m)<m-r(m)$, and $m$
and $r(m)$ are odd, therefore $r(m)\leq m-r(m)-1$,or $r(m)\leq (m-1)/2$.

 Left side of inequations follows from obvious remark that $4r^2(m)+1=m\,x < m$ is
 impossible.

$$
x(m)=\frac{4r(m)^2+1}{m}\leq
\frac{4\frac{(m-1)^2}{4}+1}{m}=\frac{m^2-2m+2}{m}=m-2+\frac{2}{m}
$$

\bigskip
{\it Lemma 5.}\label{m6} {\it $r(m)< r(m^2)$. Particularly, every number
appears in  $S$ the first time in its first degree.}

\medskip
 {\it Proof.}  Inequation   $r(m)\leq r(m^2)$ is obvious. If $r(m)=r(m^2)$ then
 $m\,x(m)=m^2 x(m^2)$ and  $x(m)=m\, x(m^2)\geq m$, which contradict to  lemma 5.

\bigskip
Put on $\nu(d)$ the number of prime factors of number  $d$.
 Let be  $d \in PP$ . Number of integer solutions of inequation
$k\,d-r_i(p)\leq n$ is
$$\left[\frac{n+r_i(d)}{d}\right]$$

 Sum $$
\sum_{i=1}^{2^{|d|}}\left[\frac{n+r_i(d)}{d}\right]
$$

 is equal to numbers divided by  $d$ among the first $n$ terms of  $S$.

Let $N=p_1 p_2\dots p_k$, so the sum
\begin{equation}\label{f1}
\sum_{d|N}\sum_{i=1}^{2^{|d|}}\left[\frac{n+r_i(d)}{d}\right] (-1)^{|d|}
\end{equation}

is equal to number of integers in sequence
 $S_1, S_2, \dots S_n$ which are mutually disjoint with $N$.

\bigskip
{\it Main Theorem} \label{m8} {\it The set of prime integers of a kind $n^2+1$
is infinite.}

\medskip
{\it Proof.} Consider the first  $k$ prime numbers  $p_1,p_2,\dots,p_k$ which
are met in sequence  $S$ as factors and put  $N=p_1 p_2\dots p_k$.  Let $n=N$
in (\ref{f1}) . In this case   $d|n$ and $r(d)<d$, therefore we can take away
all the square brackets and all items  $r_i(d)$. Then the number of integers in
$S$  which are mutuality disjoint with  $N$ is equal

 $$
 \sum_{d|N}\frac{N(-1)^{\nu(d)}2^{\nu(d)}}{d}=N\sum_{d|N}\frac{(-2)^{\nu(d)}}{d}=
 N\prod_{i=1}^k\left(1-\frac{2}{p_i}\right)=\prod_{i=1}^k(p_i-2)>0
$$

In this way after strike out the first $k$ of prime numbers $\ref{m2}$ and all
multiples  of them many of not striking rest and according to lemma 2 first of
resting  numbers is prime.

{\it Proof (second).} The number $M=4(p_1p_2\dots p_k)^2+1)$ is mutually
disjoint with $p_1, p_2, \dots p_k $ then can't be strike out from $S$.
\renewcommand{\refname}{References}
\bibliographystyle{}

\bigskip

Govorov Valentin E.

Chernomorsky bul. 4, 361.

 Moscow \quad 113556

 Russia

 e-mail:\, buch8035@mail.ru
\end{document}